\documentclass[%
twocolumn
]{mpi2015-cscpreprint}

\usepackage{graphicx,xcolor}      % include this line if your document contains figures

\usepackage{amsmath, amssymb}
\usepackage{tikz}
\usetikzlibrary{arrows,decorations.pathmorphing,backgrounds,fit,positioning,shapes.symbols,chains,mindmap,trees,automata}

\usepackage{pgfplots,amsmath}
\usepackage[absolute,overlay]{textpos}
\usepackage{color}

 \definecolor{green2}{rgb}{0, 0.5, 0}
  \definecolor{gray}{rgb}{0.5, 0.5, 0.5}

\newcommand{\bff}{{\mathbf f}}

\newcommand{\br}{{\mathbf r}}
\newcommand{\bu}{{\mathbf u}}
\newcommand{\bv}{{\mathbf v}}
\newcommand{\bw}{{\mathbf w}}
\newcommand{\bx}{{\mathbf x}}
\newcommand{\by}{{\mathbf y}}
\newcommand{\bl}{{\mathbf l}}
\renewcommand{\br}{{\mathbf r}}

\def\IR{{\mathbb R}}
\def\IC{{\mathbb C}}

\def\IL{{\mathbb L}}

\newcommand{\sIL}{{{{\mathbb L}_s}}}
\newcommand{\balpha}{{\boldsymbol{\alpha}}}

\providecommand{\Htran}[0]{\mathbf{H}} %
\providecommand{\Hreal}[0]{\mathcal{S}} %
\providecommand{\rank}[0]{\textbf{rank}} %
\providecommand{\eg}[0]{\textit{e.g.}~} \providecommand{\ie}[0]{\textit{i.e}~} %

\newenvironment{customlegend}[1][]{%
	\begingroup
	\csname pgfplots@init@cleared@structures\endcsname
	\pgfplotsset{#1}%
}{%
	\csname pgfplots@createlegend\endcsname
	\endgroup
}%

\def\addlegendimage{\csname pgfplots@addlegendimage\endcsname}

\begin{document}

	\title{On Loewner data-driven control for infinite-dimensional systems}

\author[$\ast$]{Ion Victor Gosea}
\affil[$\ast$]{Max Planck Institute for Dynamics of Complex Technical Systems, Magdeburg, Germany.\authorcr
	\email{gosea@mpi-magdeburg.mpg.de}, \orcid{0000-0003-3580-4116}}

\author[$\dagger$]{Charles Poussot-Vassal}
\affil[$\dagger$]{Information Processing and Systems Department, ONERA, Toulouse, France.\authorcr
	\email{charles.poussot-vassal@onera.fr}, \orcid{0000-0001-9106-1893}}

\author[$\ddagger$]{Athanasios C. Antoulas}
\affil[$\ddagger$]{Electrical and Computer Engineering (ECE) Department, Rice University, Houston, USA, 
	
	Max Planck Institute, Magdeburg, Germany, and Baylor College of Medicine, Houston, USA.\authorcr
	\email{aca@rice.edu}}

\shorttitle{Loewner data-driven control}
\shortauthor{I.V. Gosea, C. Poussot-Vassal, and A.C. Antoulas}
\shortdate{}


\abstract{         
	In this paper, we address extensions of the Loewner Data-Driven Control (L-DDC) methodology. First, this approach is extended by incorporating two alternative approximation methods known as Adaptive-Antoulas-Anderson (AAA) and Vector Fitting (VF). These algorithms also include least squares fitting which provides additional flexibility and enables possible adjustments for control tuning. Secondly, the standard model reference data-driven setting is extended to handle noise affecting the data and uncertainty in the closed-loop objective function. These proposed adaptations yield a more robust data-driven control design.
	%To the authors perspective, these improvements go the direction of a more robust data-driven control design. 
}
	
	\maketitle

%===============================================================================
\section{Introduction}

%{\color{red} I copied this intro to the big paper ...please feel free to modify as you consider fit.}
\subsection{From linear model-driven control to data-driven control}

In many branches of engineering, in order to satisfy accuracy requirements, the models under consideration might have large dimension and, hence, are difficult to use for control design, numerical simulations or analysis. That is why, it is of critical importance to find reliable reduced-order surrogate models instead. The latter may then be used in place of the original one. In this first case, \emph{model reduction} typically refers to a class of methodologies used for reducing the computational complexity of large-scale models of dynamical systems. The goal generally is to approximate the original model with a smaller and simpler one, having the same structure and similar response characteristics as the original. For an overview of model reduction methods, we refer the reader to the books of \cite{ACA05,BOCW17,ABG20}. Moreover, in some applications, a mathematical description of the system is not always available or involves even more complex equations. This is the case when dynamical models are described \eg by  a dedicated simulator, from which the input-output map is not available but can be evaluated. In this second case, instead of relying on equations derived from physical laws, one can infer properties and a model directly from the data, which can be done by \emph{model approximation}, see \eg \cite{ABG20}. In both cases, model reduction and approximation play the pivotal enabler role for model-driven control design.

When considering the case in which only a \emph{simulator}, or an \emph{experimental test benchmark} are accessible, instead of following the model-driven approach, one may need to use a Data-Driven Control (DDC) design rationale instead. One substantial advantage of this control tuning family is that it provides a controller tailored to the considered system, and skips the modeling phase\footnote{Note that, in many applications, the model only serves the control design and analysis, and can practically always be amended.}. Among the multiple data-driven control design approaches, we mention Virtual Reference Feedback Tuning (VRFT), introduced by \cite{Campi:2002}. The latter is particularly easy to deploy and is considered in this work. We also refer the reader to \cite{Formentin:2014}, for comparisons between model-driven and data-driven control or to \cite{Karimi:2017}, for a data-driven frequency-domain approach. Following the VRFT philosophy, belonging to the model reference methods, the control design problem is recast as an \emph{identification/data-driven approximation} one.  

%With the ever-increasing availability of measured data in many engineering fields, the need for incorporating measurements in the modeling and controlling stage has steadily grown over the last decades. The main challenge consists in using the available data effectively, so that to construct models and controllers that can accurately represent the underlying control system. In this latter case, the controller has to be designed based on experimental measurements, instead of a model description. As opposed to model-driven indirect methods, data-driven strategies directly compute the controller from experimental data. Such techniques are also known as direct methods and can be useful when such a control-oriented model is either too complex or too costly to obtain. 
\subsection{Control design via data-driven approximation}

The Loewner Data-Driven Control (L-DDC) algorithm, originally introduced in \cite{KergusWC:2017}, is a model reference technique based on frequency-domain data. Extensions of L-DDC include dealing with stability, reference model selection and controller validation (see \eg \cite{KergusLCSS:2019}). Recently, a hybrid version has also been proposed by \cite{VuilleminWC:2020}. The Loewner Framework (LF) is a data-driven model \emph{identification} and \emph{reduction technique} that was originally introduced in \cite{MA07}. Using only measured data, it directly constructs surrogate models, by employing low computational effort. For a tutorial paper on LF for linear systems, we refer the reader to \cite{ALI17}. An extension that uses time-domain data is given in \cite{PGW17}, while an extension for certain classes of nonlinear systems, \eg described by bilinear models, is given in \cite{AGI16}. The Adaptive-Antoulas-Anderson (AAA) algorithm, originally introduced in \cite{NST18}, is a data-driven rational approximation method that combines interpolation (as does LF), and also least squares fitting. AAA can also be used as a LTI modeling method since it yields a reduced-order rational function that can be interpreted as the transfer function of the surrogate reduced order model (ROM).  The AAA algorithm has recently extended for modeling of parametrized dynamics in \cite{CS20}, and for approximation of matrix-valued functions in \cite{GG20}. Finally, the Vector-Fitting (VF) method is based solely on least squares approximation, and can be also applied for surrogate modeling design.

\subsection{Paper contribution and structure}

In this paper, the L-DDC approach is for the  first time extended and compared with the VF and AAA procedures, leading to both VF-DDC and AAA-DDC control design methods. This comparison, presented by means of an infinite-dimensional irrational model representing a linear partial differential equation set, highlights the properties of each technique. We believe that this first contribution may serve practitioners in choosing an approach accordingly to the setup. Secondly, based on the usage of mixed interpolation/least squares approaches detailed in what follows, we also propose an approach to deal with the noise on the collected data and the variability of the expected performances. This is a first step toward an uncertain framework for this class of data-driven approximation DDC tuning approaches. %The essential benefit of these approaches 

The paper is organized as follows: Section \ref{sec:pb} recalls the standard DDC problem
and suggests an extension for cases with uncertainty. Section \ref{sec:ddi} then provides a complete review of the three data-driven identification methods considered here: LF, AAA and VF. Practical considerations are specifically pointed to provide as much self-contained reading as possible. Section \ref{sec:numerics} then illustrates the three algorithms proposed, first on an academic finite-order linear model, and afterwards on a  more complex model. This is ruled by linear partial differential equations, representing a transport equation phenomena. Conclusions and future research directions are discussed in Section \ref{sec:conc}.

%%%%%%%%%%%%%%%%%%%%%%%%%%%%%
\section{Problem formulation}
\label{sec:pb}

\subsection{Frequency-domain DDC (standard) rationale}

\begin{figure}[h]
\centering
\scalebox{.9}{\tikzstyle{block} = [draw, thick,fill=blue!20, rectangle, minimum height=3em, minimum width=5em,rounded corners]
\tikzstyle{sum} = [draw, thick,fill=blue!20, circle, node distance=1cm]
\tikzstyle{input} = [coordinate]
\tikzstyle{output} = [coordinate]
\tikzstyle{pinstyle} = [pin edge={to-,thick,black}]
\tikzstyle{connector} = [->,thick]

% The block diagram code is probably more verbose than necessary
\begin{tikzpicture}[auto, node distance=2cm,>=latex']
    % We start by placing the blocks
    \node [input, name=input] {};
    \node [sum, right of=input] (sum) {};
    \node [block, right of=sum] (controller) {$\mathbf K$};
    %\node [block, right of=controller, pin={[pinstyle]above:Disturbances},node distance=3cm] (system) {System};
    \node [block, right of=controller, node distance=3cm] (system) {$\Htran$};
    \node [sum, right of=system, node distance=2cm] (sum2) {};
    \node [block, above right of=controller, node distance=2cm and 2cm] (obj) {$\mathbf M$};
    % We draw an edge between the controller and system block to 
    % calculate the coordinate u. We need it to place the measurement block. 
    \draw [connector] (controller) -- node[name=u] {$\bu$} (system);
    \node [output, right of=sum2, node distance=1cm] (output) {};
    %\node [block, below of=u] (measurements) {Measurements};

    % Once the nodes are placed, connecting them is easy. 
    \draw [connector] (input) -- node {${r}$} (sum);
    \draw [connector] (sum) -- node {$e$}(controller);
    \draw [connector] (system) -- node [name=y] {$\by$}(sum2);
    %\draw [->] (y) |-| (measurements);
    %\draw [->] (measurements) -| node[pos=0.99] {$-$} node [near end] {$y_m$} (sum);
    \draw [connector] (sum2)+(-0.3cm,0) -- ++(-0.3cm,-1cm) -| node [near start] {} (sum.south);
    \draw [connector] (sum2) -- node {$\varepsilon$}(output);
    \draw [connector] (input)+(0.3cm,0) |- (obj) -| (sum2);
%    \pause
%    \node [block, above of=u, node distance=1cm,green!30,text=black,thick,draw=black,rounded corners] (model) {Reference model};
%    \draw [connector,green!30] (output)+(-0.3cm,0) |- node [near start] {} (model.east);
%    \draw [-,thick,green!30] (model.west) -| node [near start] {} (input);
    
\end{tikzpicture}}
\caption{Data-driven control problem formulation: $\mathbf{M}$ is the reference model (objective) and $\mathbf{K}$ the controller to be designed.}
\label{fig:problem_formulation}
\end{figure}
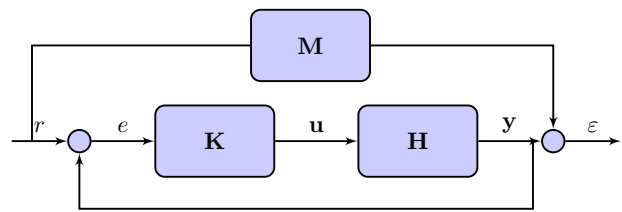

The DDC approach discussed in this paper is based on the original contribution of \cite{Campi:2002}, which was recently extended to the frequency-domain, in \cite{KergusLCSS:2019}. As depicted in Figure \ref{fig:problem_formulation}, the system to be controlled is described by $\Htran$. This latter is considered as unknown while frequency-domain input-output data are accessible such that, for pulsation $\omega_i\in\mathbb R_+$, $i=1,\dots , n$ ($n\in \mathbb N$),
\begin{equation}
    %\Htran(\imath\omega_i) = \dfrac{\tilde y(\imath\omega)}{\tilde u(\imath\omega)}, 
    \Phi_i= \dfrac{\overline \by(\imath\omega)}{\overline \bu(\imath\omega)},
\end{equation}
where $\imath=\sqrt{-1}$ and, $\overline \bu$ and $\overline \by$ refer to the Fourier transform signal of $\bu$ and $\by$, respectively. Note that $\Phi_i=\Htran(\imath\omega_i)$ holds in the ideal noise-free case. Then, function $\mathbf M$ represents the so-called objective closed-loop transfer function. This latter defines the expected response that the user needs to impose to the system when the controller is inserted in the looped architecture considered. 

The objective is to find a controller $\mathbf{K}$ minimizing the difference between the resulting closed-loop and the reference model $\mathbf{M}$. This is made possible through the definition of the so-called ideal controller $\mathbf{K}^\star$, being the LTI controller that would give the desired reference model frequency-domain behavior, if inserted in the closed-loop. This latter is defined as follows:
\begin{equation}
    \mathbf{K}^\star(\imath\omega_i)=\Phi_i^{-1}\mathbf{M}(\imath\omega_i)(I-\mathbf{M}(\imath\omega_i))^{-1}.
    \label{Kideal}
\end{equation}

%In the approach defended by \cite{KergusLCSS:2019}, authors propose to use the LF to actually compute the controller and its approximation $\mathbf K$ that fit (interpolate) the data defined in \eqref{Kideal}. One major benefit of using LF with respect to classical time-domain identification approach is that it allows finding both zeros and poles of the controller in a one shot approach.

Finding a controller $\mathbf K$ that fits $\mathbf{K}^\star(\imath\omega_i)$ can be considered to be an identification problem. In this work, this latter is considered through the lens of both data-driven \emph{interpolatory}, and of \emph{least squares} methods. This is done by comparing the methods LF, AAA, and VF, being purely interpolatory, interpolatory + least squares, and purely least squares methods, respectively (see also Section \ref{sec:ddi}).

\subsection{Frequency-domain DDC uncertain rationale}

The $\mathbf K^\star$ controller is uniquely defined by the sought-after objective closed-loop function $\mathbf M$. However, in practical applications, the following issues need to be taken into consideration. First, system's data $\Phi_i$ may be corrupted by noise, thus $\Phi_i=\Htran(\imath\omega_i)+n_i$, where $n_i\in\mathbb R$ represents the noise affecting the data\footnote{Note that in experimental setup, noise directly comes from the sensor accuracy while in the simulator-based case, it may come from the numerical Fourier transform and simulator variability.}. Second, the objective function $\mathbf M$ is not necessarily unique, and may be instead described by a set of objective functions $\{\mathbf M_j\}_{j=1}^{n_s}$ (deemed as functional for the considered process). The problem is then restated as (with $i=1,\dots,n$ and $j=1,\dots , n_s$):
\begin{equation}
    \mathbf{K}^\star(\imath\omega_i)=\Phi_i^{-1}\mathbf{M}_j(\imath\omega_i)(I-\mathbf{M}_j(\imath\omega_i))^{-1}.
    \label{Kideal_unc}
\end{equation}
Relation \eqref{Kideal_unc} will be used to address \emph{uncertainty} or \emph{robustness} issues. Here, one seeks again for controller $\mathbf K$ that fits $\mathbf{K}^\star(\imath\omega_i)$. Given this extended uncertainty problem, we next provide a review of the considered frequency-domain identification and approximation techniques: LF, AAA and VF. Each of these methods will then be embedded in Algorithm 1 of \cite{KergusLCSS:2019}, resulting in L-DDC, AAA-DDC and VF-DDC.

%%%%%%%%%%%%%%%%%%%%%%%%%%%%%
\section{Data-driven identification} 
\label{sec:ddi}

\subsection{Loewner framework (LF) interpolation}\label{subsec:Loewner}

In this section, the Loewner framework is recalled for the multi-input multi-output ({MIMO}) case. For a complete description, we refer the reader to \cite{ALI17}, and to \cite{Antoulas16} for insight in the rectangular case. Under mild considerations, the Loewner approach is a data-driven method aimed at building a rational descriptor LTI dynamical model $\Htran_m$ of dimension $m$ which interpolates given complex data, here generated by a model $\Htran$. Let the left (or row) data be given together with the right (or column) data, as below
\begin{equation}
\left.
\begin{array}{c}
(\mu_j,\bl_j^H,\bv_j^H) \\
\text{for $j=1,\dots,m$}
\end{array}
\right\}
\text{~~and~~}
\left\{
\begin{array}{c}
(\lambda_i,\br_i,\bw_i) \\
\text{for $i=1,\dots ,m$}
\end{array}
\right. ,
\label{eq:loewnerInput}
\end{equation}
where $\bv_j^H=\bl_j^H\Htran(\mu_j)$ and $\bw_i=\Htran(\lambda_i)\br_i$, with $\bl_j\in\IC^{n_y\times 1}$, $\br_i\in\IC^{n_u\times 1}$, $\bv_j\in\IC^{n_u\times 1}$ and $\bw_i\in\IC^{n_y\times 1}$. In addition, the set of distinct interpolation points $\{ z_k\}_{k=1}^{2m} \subset \IC$ is split up into two equal subsets ($\lambda_i, \ \mu_j\in\IC$), \ie
\begin{equation}\label{eq:shift}
\{z_k\}_{k=1}^{2m} =\{\mu_j\}_{j=1}^{m} \cup \{\lambda_i\}_{i=1}^{m}.
\end{equation}
The method then consists in building the \emph{Loewner} matrix $\IL \in \IC^{m\times m}$ and \emph{shifted Loewner} matrix $\sIL \in \IC^{m\times m}$ defined as follows, for $i=1,\dots,m$ and $j=1,\dots,m$:
\begin{align}\label{eq:loewnerMatrices}
\begin{split}
[\IL]_{j,i} &= \dfrac{\bv_j^H\br_i - \bl_j^H\bw_i}{\mu_j - \lambda_i} 
= \dfrac{\bl_j^H\big( \Htran(\mu_j) - \Htran(\lambda_i) \big) \br_i}{\mu_j - \lambda_i}, \\
\,[\sIL]_{j,i} &= \dfrac{\mu_j\bv_j^H\br_i - \lambda_i\bl_j^H\bw_i}{\mu_j - \lambda_i}
= \dfrac{ \bl_j^H\big( \mu_j\Htran(\mu_j) - \lambda_i\Htran(\lambda_i) \big) \br_i}{\mu_j - \lambda_i}.
\end{split}
\end{align}
Then, the model $\Htran_m$ given by the descriptor realization,
\begin{equation}
\Hreal_m:\left \lbrace
\begin{array}{rcl}
E_m  \delta\left\{\bx(\cdot)\right\} &=& A_m \bx(\cdot) + B_m \bu(\cdot)\\
\by(\cdot) &=&C_m \bx(\cdot)
\end{array}
\right. ,
\label{eq:loewnerDescrR}
\end{equation}
where $E_m = -\IL$, $A_m = -\sIL$, $[B_m]_k = \bv_k^H$ and $[C_m]_k = \bw_k$ (for $k=1,\ldots,m)$, with the related transfer function 
\begin{equation}
\Htran_m(\xi) = C_m(\xi E_m-A_m)^{-1}B^m,
\label{eq:loewnerDescrC}
\end{equation}
interpolates $\Htran$ at the given driving frequencies and directions defined in \eqref{eq:loewnerInput}, \ie satisfies the conditions
\begin{equation}
\begin{array}{rcl}
\bl_j^H\Htran_m(\mu_j) &=& \bl_j^H\Htran(\mu_j) \\% \text{ and }
\Htran_m(\lambda_i)\br_i &= &\Htran(\lambda_i) \br_i
\end{array}.
\label{eq:loewnerIntep}
\end{equation}
%\begin{remark}[About $\delta\left\{\cdot\right\}$ and $\zeta$ notations]
Note that ``$(\cdot)$'' denotes the time-domain variable considered in \eqref{eq:loewnerDescrR}; this can either be ``$(t)$'' for continuous-time models ($t\in\IR_+$), or ``$[q]$'' for discrete-time models ($q\in\mathbb Z$). Similarly, in \eqref{eq:loewnerDescrR} ``$\delta\left\{\cdot\right\}$'' stands as the shift operator being either $\delta\{\bx(t)\}=\mathbf{\dot x}(t)$ in the continuous-time case, and $\delta\{\bx(q)\}=\bx[q+1]$ in the discrete-time one. Note also that in \eqref{eq:loewnerDescrC}, $\xi$ represents the associated Laplace complex variable $\xi=s$ in the continuous-time case, and the forward shift $\xi=z$ in the discrete-time one.
%\end{remark}

Assuming that the number $2m$ of available data is large enough, then  it was shown in \cite{MA07} that a minimal model $\Htran_n$ of dimension $n < m$ (that still interpolates the data) can be computed with a projection of \eqref{eq:loewnerDescrR} provided that the following holds for $k=1,\ldots,2m$
\begin{equation}
%\scriptsize
    \rank (z_k \IL - \sIL) = 
    \rank ([\IL,\sIL]) = 
    \rank ([\IL^H,\sIL^H]^H) = n,
    \label{eq:rankCond}
\end{equation}
where $z_k$ are as in \eqref{eq:shift}. In that case, let $Y \in \IC^{m \times n}$ be the matrix containing the first $n$ left singular vectors of $[\IL,\sIL]$ and $X \in \IC^{m \times n}$  the matrix containing the first $n$ right singular vectors of $[\IL^H,\sIL^H]^H$. Then,
\begin{align}    \label{eq:proj}
\begin{split}
    E_n &= Y^H E_m X,\,
    A_n = Y^H A_m X,\\
    B_n &= Y^H B_m,\,
    C_n = C_m X,
\end{split}
\end{align}
is a realization of the model $\Htran^n$, given as, 
\begin{equation}
\Htran_n(\xi) = C_n( \xi E_n-A_n)^{-1}B_n,
\label{eq:loewnerDescrCn}
\end{equation}
with the same structure as \eqref{eq:loewnerDescrC}, encoding a \emph{minimal McMillan degree} equal to $\rank(\IL)$. The quadruple given by $\Hreal_n:(E_n,A_n,B_n,C_n,0)$ is a descriptor realization of $\Htran_n$. Note that if $n$ in \eqref{eq:rankCond} is greater than $\rank(\IL)$, then $\Htran_n$ can either have a direct-feedthrough term or a polynomial part. Finally, the number $n$ of singular vectors composing $Y$ and $X$ used to project the system $\Htran_n$ in \eqref{eq:proj} may be decreased to $r<n$ at the cost of imposing an approximate interpolation of data, leading to the reduced model $r$-th order rational model.% denoted $\Htranr(\xi)=\hat C(\xi\hat E-\hat A)^{-1}\hat B$. 
This allows a trade-off between complexity of the resulting model and accuracy of the interpolation.% (see \cite{MA07}).

\subsection{The AAA algorithm}\label{subsec:AAA}

The AAA algorithm, originally proposed in \cite{NST18}, represents an adaptive extension of the interpolation-based method introduced in \cite{AA86}. It is a robust, fast and effective  method that was mainly used for scalar rational interpolation applications. AAA is a multi-step algorithm, that computes at step $\ell$ a rational approximant of order $(\ell,\ell)$ in barycentric representation. In this note we discuss a slightly modified version from that in \cite{NST18}, in the sense that the approximant at step $\ell$ is strictly proper, \ie of order $(\ell-1,\ell)$. Additionally, as for the Loewner method, we will enforce real-valued models. Finally, we restrict the presentation to the SISO case (the MIMO case was addressed in \cite{GG20}).

As in (\ref{eq:shift}), we consider at step $\ell \geq 1$ the data splitting:
\begin{align}\label{eq:split_AAA}
\begin{split}
\text{data points}: \ \{z_k\}_{k=1}^{2m} &=\{\nu_j\}_{j=1}^{\ell} \cup \{\eta_i\}_{i=1}^{2m-\ell}, \\
\text{data values}: \ \{f_k\}_{k=1}^{2m} &= \{h_j\}_{j=1}^{\ell} \cup \{g_i\}_{i=1}^{2m-\ell}.
\end{split}
\end{align}
Note that in the representation given in (\ref{eq:split_AAA}), the values $f_k$ represent the measurements evaluated at the points $z_k$, while $h_j$ and $g_i$ are the ones evaluated at $\nu_j$, and respectively at $\eta_i$.
%$\{H(z_k)\}_{k=1}^{2m} &= \{H(\nu_j)\}_{j=1}^{\ell} \cup \{H(\eta_i)\}_{i=1}^{2m-\ell},$
The rational interpolant $\Htran_\ell$, obtained after $\ell$ iterations of the AAA algorithm, has the form
\begin{equation}\label{eq:bary}
	\Htran_\ell(\xi) =  \frac{\sum_{j=1}^\ell \frac{\alpha_j^{(\ell)} h_j}{\xi - \nu_j}}{1+\sum_{j=1}^\ell \frac{\alpha_j^{(\ell)}}{\xi - \nu_j}},
\end{equation}
with nonzero barycentric weights $\alpha_j ^{(\ell)} \in\mathbb{C}$, pairwise distinct support points $\nu_j\in\mathbb{C}$, and function values $h_j$. Based on the representation in (\ref{eq:bary}), interpolation is enforced at the first subset of data points $\{\nu_j\}_{j=1}^{\ell}$, \ie $\Htran_{\ell}(\nu_j) = h_j$ for $1 \leq j \leq \ell$. In order to completely determine the approximant $\Htran_\ell$, one needs to also find the barycentric weights $\alpha_1^{(\ell)} ,\ldots,\alpha_\ell^{(\ell)}$. This is done by solving  a least squares problem. Finally, the next support point is chosen by means of a greedy selection.

Let $\tau$ be the desired tolerance for data approximation and let $n$ denote the target dimension.
The modified AAA algorithm can be summarized as follows: 
\medskip
\begin{enumerate}
\item \textit{Initialization step}\\ Set $\ell=0$, $\Omega^{(0)} := \Omega$, and $\Htran_{0}(\xi) = \frac{1}{2m} \sum_{k=1}^{2m} f_k$.

\item \textbf{While} $\underset{1 \leq k \leq 2m}{\max} \Big{|}  f_k - \Htran_{\ell}(z_k)  \Big{|} > \tau$ \ and \ $\ell < n$ \\ 
%(\ie the relative deviation is not small enough or the order $n$ is not reached yet; if any of this criterion is met, return $H_{\ell}$ and stop).

\item \textbf{do} $\ell = \ell+1$.

\item Find $\nu_{\ell} \in \Omega^{(\ell-1)}$ so that $\nu_{\ell} = \underset{1 \leq k \leq 2m}{\text{argmax}} |f_k - \Htran_{\ell-1}(z_k)|$, with $\Htran_{\ell}(s)$ as in (\ref{eq:bary}). Set $h^\ell := \Htran(\nu_\ell)$, and also  $\Omega^{(\ell)} := \Omega^{(\ell-1)}\setminus \{\nu_\ell\}$.
\item Compute weights $\alpha_1^{(\ell)},\ldots,\alpha_\ell^{(\ell)}$ to minimize the deviation in the measurements, \ie solve the problem
\begin{equation}\label{eq:min_prob_nonlin}
    \underset{\alpha_1^{(\ell)},\ldots,\alpha_\ell^{(\ell)}}{\min}  \sum_{k=1}^{2m} ( \Htran_\ell(z_k)-f_k )^2.
\end{equation}

%\begin{align}\label{eq:min_prob_nonlin}
%\begin{split}
%\min \sum_{k=1}^{2m} \left( \frac{\sum_{j=1}^\ell \frac{\alpha_j^{(\ell)} h_j}{z_k - \nu_j}}{1+\sum_{j=1}^\ell \frac{\alpha_j^{(\ell)}}{z_k - \nu_j}} -f_k \right)^2 \\
%\Leftrightarrow \min \sum_{k=1}^{2m} \left( \frac{\sum_{j=1}^\ell \frac{\alpha_j^{(\ell)} (f_k-h_j)}{z_k - \nu_j}+f_k}{1+\sum_{j=1}^\ell \frac{\alpha_j^{(\ell)}}{z_k - \nu_j}}  \right)^2.
%\end{split}
%\end{align}
\item Instead of solving the nonlinear problem in (\ref{eq:min_prob_nonlin}), one solves a linearized problem by substituting $\Htran_\ell$ in (\ref{eq:bary}):
\begin{align}\label{eq:min_prob_lin}
\hspace{-2mm}   \underset{\balpha^{(\ell)}}{\min} \sum_{k=1}^{2m} \sum_{j=1}^\ell  \left( \frac{(f_k-h_j)\alpha_j^{(\ell)}}{z_k - \nu_j}+f_k  \right)^2 \Leftrightarrow \Vert \IL \balpha^{(\ell)} + \bff \Vert_2^2,
\end{align}
where $\IL \in \mathbb{C}^{2m \times \ell}$ with $\IL_{k,j} = \frac{f_k-h_j}{z_k - \nu_j}$ is a Loewner matrix, while $\balpha^{(\ell)} \in \mathbb{C}^{\ell}$, and $\bff \in \mathbb{C}^{2m }$.
\item Compute the solution to (\ref{eq:min_prob_lin}) as $\balpha^{(\ell)} = - \IL^{\#} \bff$, where $\IL^{\#} \in \mathbb{C}^{\ell \times 2m}$ is the pseudo-inverse of matrix $\IL$. \\
%and form the approximant $H^\ell$ as in (\ref{eq:bary}).\\
\textbf{end}
%\[
%	H(s)+\sum_{j=1}^\ell \frac{\alpha_j^{(\ell)}}{s - \nu_j} H(s) \approx \sum_{j=1}^\ell \frac{\alpha_j^{(\ell)} H(\nu_j)}{s - \nu_j}
%\]
\end{enumerate}
\medskip

It is to be noted that a realization $\Hreal_n:(I_n,A_n,B_n,C_n,0)$ of an order $(n-1,n)$ of the AAA reduced-order model  can be expressed as follows
	\begin{align*}
		{A_n} &= \text{diag}(\nu_1,\ldots,\nu_n) -  {B_n} e_n^T, \\
		 {B_n} &= \left[ \begin{array}{cccc} \alpha_1^{(n)}  & \ldots & \alpha_n^{(n)}  \end{array} \right]^T, \ \ {C_n} = \left[ \begin{array}{cccc} h_1 & \ldots & h_n \end{array} \right],
	\end{align*}
	where $e_n = \left[ \begin{matrix} 1 & \cdots & 1 \end{matrix} \right]^T$ and ${A_n}$ is rank-1 perturbation of a diagonal matrix (composed of the chosen support points). 

%	Then, the following holds:
%	\begin{equation}\label{eq:bary2}
%	H_n(s) =  \frac{\sum_{j=1}^n \frac{\alpha_j^{(n)} h_j}{s - \nu_j}}{1+\sum_{j=1}^n \frac{\alpha_j^{(n)}}{s - \nu_j}} = \hat{C} (sI_n-\hat{A})^{-1} \hat{B},
%\end{equation}

%In order to enforce a real-valued LTI model, one needs to modify the previously presented algorithm by choosing not only $\nu_\ell$ as support point, but also its conjugate $\nu^*_\ell$.

%MATLAB implementations of the original AAA algorithm can be found in \cite{NST18} and the Chebfun package~\cite{chebfun}.

\subsection{Vector fitting (VF)}

Vector fitting, originally introduced in \cite{GS99}, is an effective approximation method used for constructing rational approximants designed to fit given frequency response measurements. The method is based on least squares approximation of the data values by a rational function, using an iterative reallocation of the approximants' poles. VF computes a rational function in pole-residue format given by
% that approximates the original scalar function $f: \mathbb{C} \rightarrow \mathbb{C}$, given its samples $\{f(\lambda_1), \ldots, f(\lambda_p)\}$:
\vspace{-2mm}
\begin{equation*}
%\tilde{H}_n
\Htran_n(\xi) =  \sum_{k=1}^{n} \frac{\gamma_k}{\xi-\zeta_k}= \frac{N(\xi)}{D(\xi)}.
\vspace{-2mm}
\end{equation*}
The approximation problem is formulated as follows
\begin{equation}\label{vf_prob_nonlin}
    \underset{\gamma_k,\zeta_k}{\min}  \sum_{k=1}^{2m} ( \Htran_n(z_k)-f_k )^2 \Rightarrow     \underset{\gamma_k,\zeta_k}{\min}  \sum_{k=1}^{2m} \left( \frac{N(z_k)}{D(z_k)}-f_k \right)^2.
\end{equation}
Since the poles $\zeta_k$ enter nonlinearly in (\ref{vf_prob_nonlin}), this again represents a nonlinear problem. Instead of solving this, one introduces an iterative algorithm that is initiated by choosing the degree $n$ of the rational approximant and an initial guess for the poles $\{\zeta_1^{(0)}, \ldots, \zeta_n^{(0)} \}$. At iteration step $j \geq 0$, the goal is to determine the parameters 
$c_i^{(j)}$ and $d_i^{(j)}$ that solve the linearized problem (the poles are excluded from the variable set) 

\begin{align*}
& \underset{c_i^{(j)},d_i^{(j)}}{\min} \sum_{k=1}^{2m} (N^{(j)}(z_k) - D^{(j)}(z_k) f_k)^2 \\
& \Leftrightarrow \underset{c_i^{(j)},d_i^{(j)}}{\min} \bigg{(} \underbrace{\sum_{i=1}^{n} \frac{c_i^{(j)}}{z_k-\zeta_i^{(j-1)}}}_{N^{(j)}(\xi)} - \big( \underbrace{1+\sum_{i=1}^{n} \frac{d_i^{(j)}}{z_k-\zeta_i^{(j-1)}}}_{D^{(j)}(\xi)} \big) f_k \bigg{)}^2.
\end{align*}
%for all $1\leq k \leq 2m$.
%is solved with least squares error over  all $z\in\Lambda$. 
The problem formulated above is linear and can be hence solved directly. Afterwards, the next set of poles, given by $\{\zeta_1^{(j)}, \ldots, \zeta_d^{(j)} \}$, is computed as the roots of the numerator of $D^{(j)}(\xi)$ by solving a linear eigenvalue problem. The iteration continues until a convergence criterion is satisfied (the poles are the same up to a tolerance value). This procedure naturally leads to a realization  $\Hreal_n:(I_n,A_n,B_n,C_n,D_n)$. It is to be noted that VF is a \emph{non-interpolatory} method. Additionally, note also that the degree of the approximant computed with VF is fixed, while for AAA it increases with each iteration.

%%%%%%%%%%%%%%%%%%%%
\section{Numerical examples} \label{sec:numerics}

By following the DDC setup presented in Section \ref{sec:pb}, one naturally extends L-DDC to AAA-DDC and VF-DDC algorithms. In this section, these procedures are applied for two numerical use cases to construct a controller, as on Figure \ref{fig:problem_formulation}, that tracks some closed-loop performances. First, a simple rational model is considered in Section \ref{ssec-ex_toy}, and second, in Section \ref{ssec-ex_transport}, a more challenging irrational infinite-dimensional model is involved.

\subsection{Academic example}
\label{ssec-ex_toy}

In this first academic example, we consider a continuous-time rational LTI model $\Htran$ described by the following realization $(E,A,B,C,D)=(1,-1,0.5,1,0)$. The input-output open-loop data $\Phi_i$ are collected for $n=60$ pulsations points $\omega_i$ sampled from $10^{-2}$ to $10^2$ with a logarithmic spacing. The objective function is first set to $\mathbf M(s)=\frac{1}{s^2/p^2+2s/p+1}$ (with $p=1$). By considering the "standard" case  in \cite{KergusLCSS:2019}, one seeks a controller that fits the ideal one defined by \eqref{Kideal}. Based on the rank conditions given in Section \ref{subsec:Loewner}, the Loewner procedure indicates that a controller with order 2 is sufficient to match the behavior. By applying the three procedures in Section \ref{sec:ddi}, the following 2-nd order controllers are computed:%, for L-DDC, AAA-DDC and VF-DDC:
\vspace{-2mm}
\begin{equation}
\begin{array}{rcl}
    \mathbf K_{\text{loe}}(s) &=& \dfrac{2 s + 2}{s^2 + 2 s - 4.441\cdot 10^{-16}},\\
    \mathbf K_{\text{aaa}}(s) &=& \dfrac{2 s + 2}{s^2 + 2 s - 2.602\cdot 10^{-16}},\\
    \mathbf K_{\text{vf}}(s) &=& \dfrac{2 s + 2}{s^2 + 2 s + 1.327\cdot 10^{-16}},
\end{array}
\label{eq:toy_classic}
\end{equation}
which all ensure stable closed-loops and exactly lead to the performances dictated by $\mathbf M$.%, as illustrated on Figure \ref{fig:toy1}.
%\begin{figure}
%    \centering
%    \includegraphics[width=\columnwidth]{toy_classic.pdf}
%    \caption{Classic case (from top left to bottom right): Bode gain, step response, Bode phase and controller gain.}
%    \label{fig:toy1}
%\end{figure}

\begin{figure*}[h]
\begin{tikzpicture}
\begin{customlegend}[legend columns=-1, legend style={/tikz/every even column/.append style={column sep=0.4cm}} , legend entries={
	%\footnotesize{BT-allswitch},
	 \normalsize{Reference model},\normalsize{Open-loop}, \normalsize{VF-DDC}, \normalsize{AAA-DDC}, \normalsize{L-DDC}}]
\addlegendimage{color=gray,solid,line width=2.5pt,forget plot}
\addlegendimage{color=blue,solid,line width=2.5pt,forget plot}
\addlegendimage{color=cyan,dotted,line width=2.5pt,forget plot}
\addlegendimage{color=green2,solid,line width=2.5pt,forget plot}
\addlegendimage{color=red,dashed,line width=2.5pt,forget plot}
\end{customlegend}
\end{tikzpicture}
\centering
    \includegraphics[width=1.5\columnwidth]{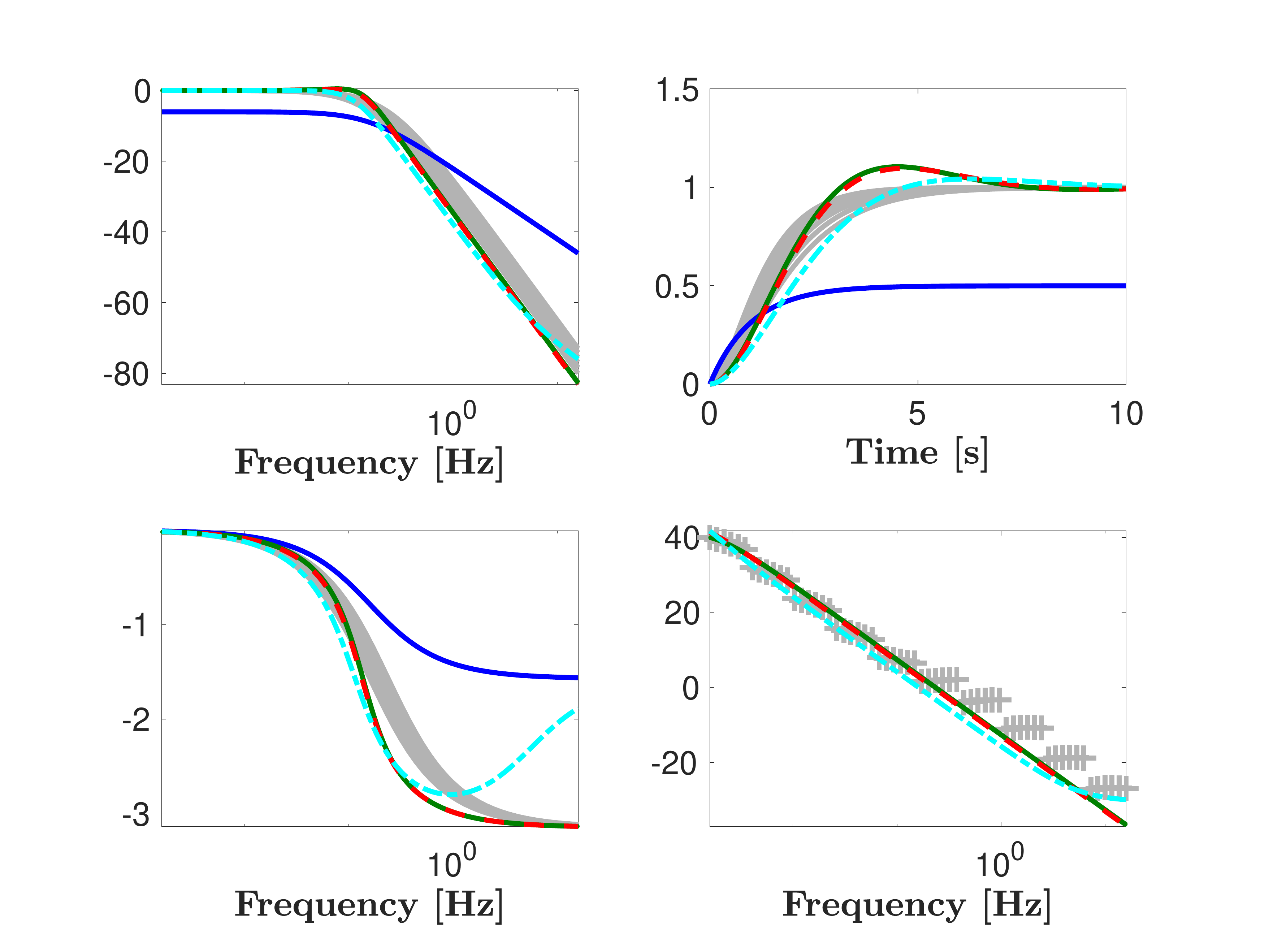}
        \vspace{-4mm}
    \caption{Uncertain case (from top left to bottom right): Bode gain, step response, Bode phase and controller gain.}
    \label{fig:toy2}
\end{figure*}

Next, instead of using an unique function $\mathbf M(s)$, one uses a family of objective functions with similar form given by $\mathbf M_j(s)$, where $j=1,\dots ,n_s$ ($n_s=6$) and $p$ is varied as $p=[1,1.1,1.2,1.3,1.4,1.5]$. This corresponds to the so-called "uncertain" case treated in \eqref{Kideal_unc}. Then, as the number of samples $n$ of $\Phi_i$ remains constant, the pulsation grid is split in $n_s$ sub-grids (the way to systematically subdivide is still an open question). Applying the three control algorithms leads then to the results given in Figure \ref{fig:toy2}, obtained with 2-nd order controllers given as:
\vspace{-1mm}
\begin{equation}
\begin{array}{rcl}
    \mathbf K_{\text{loe}}(s) &=& \dfrac{1.405 s + 0.0002368}{s^2 + 0.0003097 s - 2.114\cdot 10^{-5}},\\[3mm]
    \mathbf K_{\text{aaa}}(s) &=& \dfrac{1.469 s + 0.0004291}{s^2 + 0.0005025 s - 4.687\cdot 10^{-5}},\\[3mm]
    \mathbf K_{\text{vf}}(s) &=& \dfrac{0.03021 s^2 + 1.012 s + 0.01233}{s^2 + 0.01114 s + 3.072\cdot 10^{-5}}.
\end{array}
\label{eq:toy_uncertain}
\end{equation}

In the "standard" case, when one single objective is sought, all methods provide a similar controller. However, the Loewner approach is the only one able to a priori detect the correct order. In the "uncertain" case, when the family of objective behaviors $\mathbf M_j$ is set, all three methods show to perform well (see Figure \ref{fig:toy2}). Interestingly, L-DDC and AAA-DDC are very similar. This can be justified by the choice of the AAA interpolation points, located in the high-gain range (low frequencies) of the controller objective. Indeed, in this specific case, as the  controller integral action represents the most important energy, it is not surprising that both L-DDC and AAA-DDC focus mostly on it. 
%(this is also a direction for future developments).
As expected, the VF-DDC based on least squares fit provides a compromise. This aspect will further be discussed in Section \ref{ssec-ex_transport}. 
%In the next section, a more challenging example is used.

%This order is then used in the AAA-DDC and VF-DDC, leading to similar results. 

\vspace{-1mm}

\subsection{Transport phenomena use-case}
\label{ssec-ex_transport}

\subsubsection{Original problem description}
\hspace{-3mm} The second example involves a one-dimensional transport equation controlled at its left boundary. This phenomena is represented by a linear PDE with constant coefficients, as described in \eqref{eq:exTransportTime}. 
\begin{equation}
\begin{array}{rcll}
&& \hspace{-9mm} \dfrac{\partial \tilde y(x,t)}{\partial x} + 2x\dfrac{\partial \tilde y(x,t)}{\partial t} = 0 & \text{~~(transport equation)}\\ 
&& \hspace{-9mm} \tilde y(x,0)=0& \text{~~(initial condition)} \\
&&\hspace{-9mm} \tilde y(0,t)=\dfrac{1}{\sqrt{t}}*\tilde u_f(0,t)& \text{~~(boundary control)}\\ 
&& \hspace{-9mm} \dfrac{\omega_0^2}{s^2+m\omega_0s+\omega_0^2}u(0,s)=u_f(0,s) &\text{~~(actuator model)},
\end{array}
\label{eq:exTransportTime}
\end{equation}
where $x\in[0~L]$ ($L=3$) is the space variable, $t$ the time variable and $s$, the Laplace one. Then, $\omega_0=3$ and $m=0.5$ are the input actuator parameters. The scalar input of the model is the vertical force applied at the left boundary, \ie at $x=0$. We denote the input $\tilde u(0,t)$ in the time domain or $u(0,s)$ in the Laplace domain. Similarly, the output at location $x$ is given as $\tilde y(x,t)$ for the time domain and $y(x,s)$ in Laplace domain. Such a transport equation may be used to represent a simplified one-dimensional wave equation used in telecommunications, traffic jam prevention, \emph{etc}.

\subsubsection{Equivalent irrational transfer function}

By applying the Laplace transform to the transport equation, one obtains 
\begin{equation}
\dfrac{\partial  y(x,s)}{\partial x} + 2x \left(s   y(x,s) - \tilde y(x,0)\right) = 0,
\end{equation}
for which the solution can be given in closed-form by
$$
y(x,s) = a(s)e^{\int -2xs dx} = a(s)e^{-x^2s}.
$$
The boundary condition $\tilde  y(0,t)=\dfrac{1}{\sqrt{t}}*\tilde u_f(t)$ is transformed into $y(0,s)=\dfrac{\sqrt{\pi}}{\sqrt{s}} u_f(s)$, and hence we have that $a(s)=\dfrac{\sqrt{\pi}}{\sqrt{s}} u_f(s)$. The transfer function from input $u(0,s)$ to  output $y(x,s)$ reads
\begin{equation}
\begin{array}{rcl}
y(x,s) &=& \dfrac{\sqrt{\pi}}{\sqrt{s}}e^{-x^2s} \dfrac{\omega_0^2}{s^2+m\omega_0s+\omega_0^2}u(0,s) \\[3mm]
&=& \mathbf G(x,s) u(0,s).
\end{array}
\label{eq:exTransportFreq}
\end{equation}

Relation \eqref{eq:exTransportFreq} links the (left boundary) input to the output through an irrational transfer function $\mathbf G(x,s)$ for any $x$ value\footnote{Interestingly, the exact time-domain solution of \eqref{eq:exTransportTime}, along $x$, is given by $\tilde y(x,t)=\tilde u_f^{t-x^2}/\sqrt{t}$, where $\tilde u_f$ is the output of the second order actuator transfer function, in response to $u$.}. For illustration purpose, let us now consider that one single sensor is available, and is located at $x_m=1.9592$ along the $x$-axis\footnote{In the rest of the paper, $x$ will be discretized with 50 points from 0 to $L=3$, and $x_m$ has been chosen to be located at $x(\lfloor{ 50 \times 2/3}\rfloor )$.}. The transfer from the same input $u(0,s)$ to $y_{x_m}(s)=y(x_m,s)$ is then given by
\begin{equation}
y_{x_m}(s)  =\mathbf H(s) u(0,s),
\label{eq:exTransportFreq_xm}
\end{equation}
where $\mathbf H$ is now a one output one input  transfer function. %complex-valued

\subsubsection{Control objective and design}

\begin{figure*}[h]
\begin{tikzpicture}
\begin{customlegend}[legend columns=-1, legend style={/tikz/every even column/.append style={column sep=0.4cm}} , legend entries={
	%\footnotesize{BT-allswitch},
	 \normalsize{Reference model},\normalsize{Open-loop}, \normalsize{VF-DDC}, \normalsize{AAA-DDC}, \normalsize{L-DDC}}]
\addlegendimage{color=gray,solid,line width=2.5pt,forget plot}
\addlegendimage{color=blue,solid,line width=2.5pt,forget plot}
\addlegendimage{color=cyan,dotted,line width=2.5pt,forget plot}
\addlegendimage{color=green2,solid,line width=2.5pt,forget plot}
\addlegendimage{color=red,dashed,line width=2.5pt,forget plot}
\end{customlegend}
\end{tikzpicture}
    \centering
    \includegraphics[width=1.5\columnwidth]{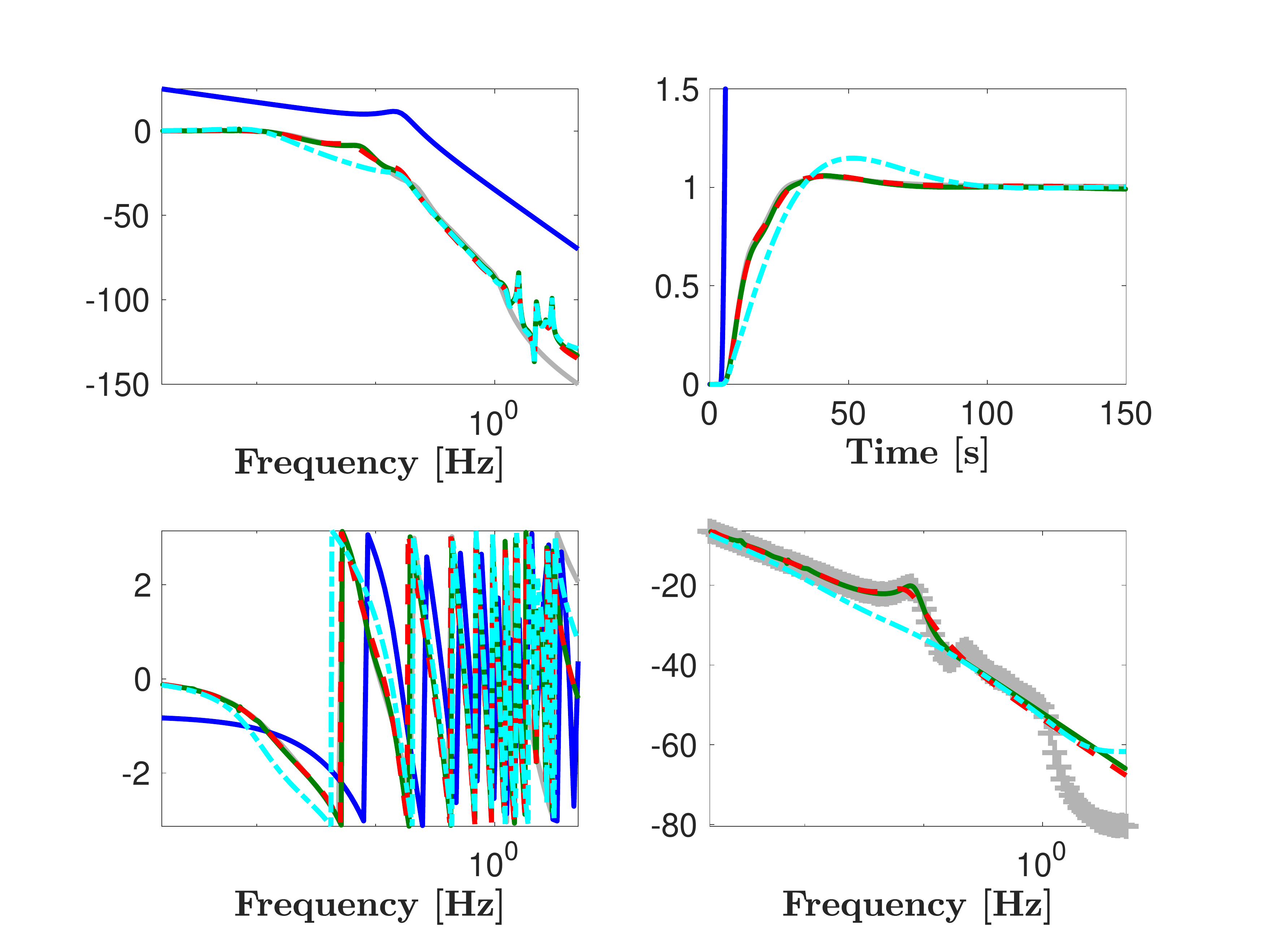}
    \vspace{-4mm}
    \caption{Standard case (from top left to bottom right): Bode gain, step response, Bode phase and controller gain.}
    \label{fig:notAtoy1}
\end{figure*}

The transport phenomena of $\Htran$ is irrational, delayed and has a limit of stability singularity. The objective of the control is to stabilize and provide some closed-loop performances. The considered measurements $\Phi$ are computed from $\Htran(\imath\omega_i)$, for $n=100$ pulsations $\omega_i$ collected between $10^{-2}$ and  $10^{1.5}$ with a logarithmic spacing. Following the control architecture of Figure \ref{fig:problem_formulation}, the data-driven control methods presented are now evaluated. Without entering into details, due to the system physical limitations, the considered reference model $\mathbf M$ is an input delayed model with oscillatory behavior\footnote{The $\mathbf M$ transfer, together with the code will be provided in the final version of the paper.}, filtered with a first order model $1/(s/p+1)$ with different parameters $p=0.1$ (in the "standard" case) and $p=p_j$ (for 5 linearly-spaced $p_j$ between $0.05$ and $0.2$) in the "uncertain" one. Similarly to the previous case, considering the "standard" problem, the L-DDC methods indicates that a 14-th order controller is enough to achieve the desired performances. Inserting the controller in the closed-loop then leads to the results depicted in Figure \ref{fig:notAtoy1}. Here again, L-DDC provides the exact expected performances.

\begin{figure*}[h]
\begin{tikzpicture}
\begin{customlegend}[legend columns=-1, legend style={/tikz/every even column/.append style={column sep=0.4cm}} , legend entries={
	%\footnotesize{BT-allswitch},
	 \normalsize{Reference model},\normalsize{Open-loop}, \normalsize{VF-DDC}, \normalsize{AAA-DDC}, \normalsize{L-DDC}}]
\addlegendimage{color=gray,solid,line width=2.5pt,forget plot}
\addlegendimage{color=blue,solid,line width=2.5pt,forget plot}
\addlegendimage{color=cyan,dotted,line width=2.5pt,forget plot}
\addlegendimage{color=green2,solid,line width=2.5pt,forget plot}
\addlegendimage{color=red,dashed,line width=2.5pt,forget plot}
\end{customlegend}
\end{tikzpicture}
   \centering
    \includegraphics[width=1.5\columnwidth]{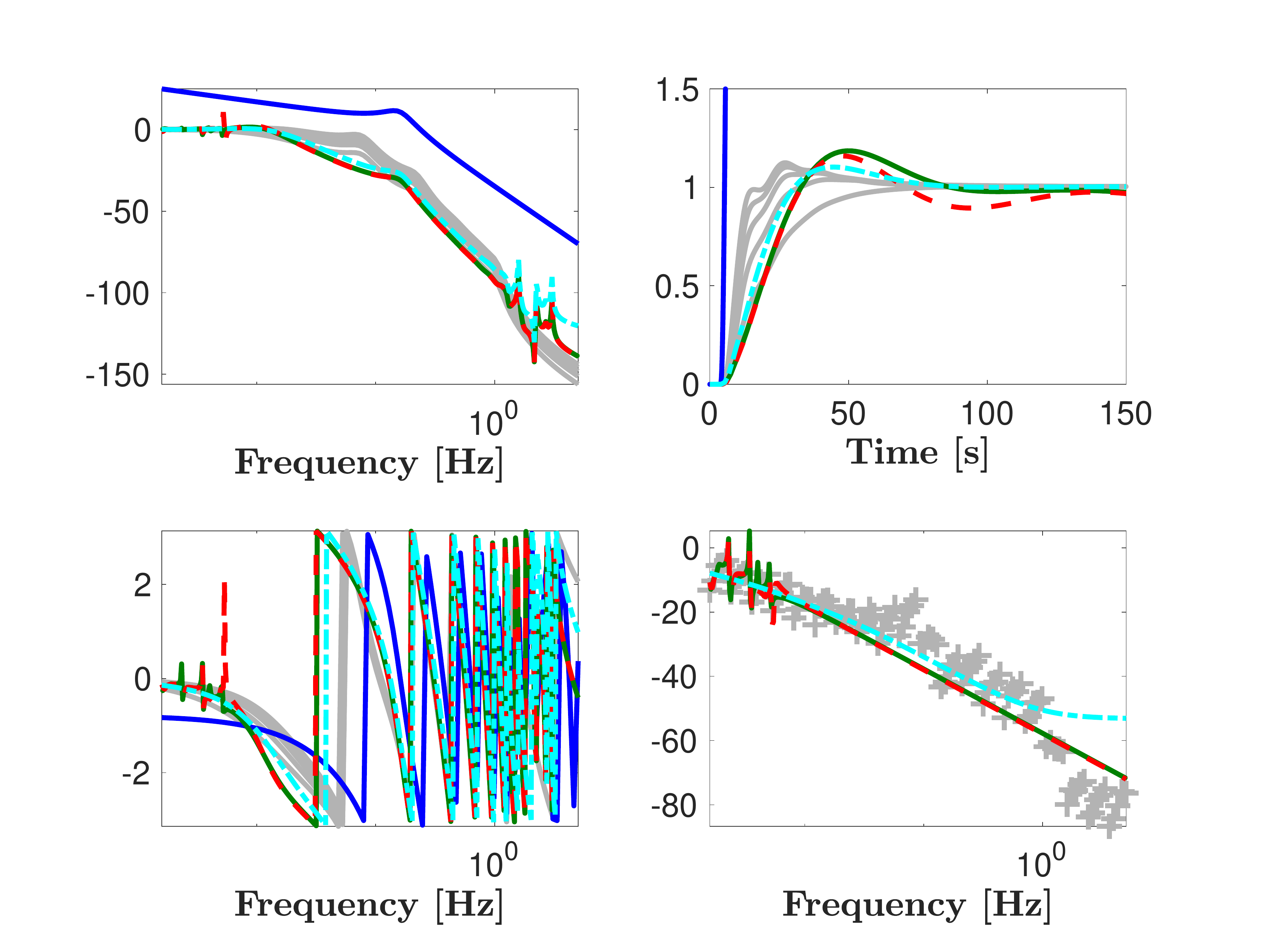}
        \vspace{-4mm}
    \caption{Uncertain case (from top left to bottom right): Bode gain, step response, Bode phase and controller gain.}
    \label{fig:notAtoy2}
\end{figure*}

As rooted on this first result, the "uncertain" case is now treated and leads to results presented in Figure \ref{fig:notAtoy2}. In this second case, noise is also added on the collected data, considering $\Phi_i (1+n_i)$ instead of $\Phi_i$ (where $n_i$ is a randomly generated number between 0 and 0.5).

In both the "standard" and "uncertain" cases, the system is stabilized, using only data. This is achieved while avoiding the modeling step and/or a dedicated work on the PDE simulator. Each of the L-DDC, AAA-DDC and VF-DDC  methods provide satisfactory performances. The Loewner-driven one has a considerable practical advantage, by providing the controller order. Moreover, it shows also to provide a good compromise between all $\mathbf M_j$. 
%However, to the authors experience, a balanced comment should be added:
Here, the VF-DDC appears to be more robust when addressing a family of objective functions $\mathbf M_j$, which is not surprising. AAA-DDC represents an interesting trade-off between the two approaches as it blends interpolation-based and least squares methods. 

\vspace{-2mm}

%%%%%%%%%%%%
\section{Conclusion}\label{sec:conc}
In this paper, the frequency-domain L-DDC rationale is revisited, first with two additional identification methods, leading to the AAA-DDC and VF-DDC algorithms. These two algorithms provide the user with alternative solutions, representing a trade-off between interpolation and least squares approximation. In addition, the original identification problem is extended to handle a family of objective functions. This allows dealing both with robustness issues, and with allowing additional degrees of freedom. Due to space limitations, a complete comparison is not set here, but mostly pointed out to. This comparison and identification methods adaptation will be addressed in future works. As for most DDC methods, stability assessment still remains an open issue. In addition, the reference model selection is still a topic for further research. Future works will address this.% following \eg, the interpolatory-based approach proposed in \cite{PoussotGTA:2020}. 

\small

\bibliographystyle{spmpsci}  

\bibliography{GPA_ECC_2021_ref}             % bib file to produce the bibliography

\end{document}